\documentclass[preprint,12pt]{elsarticle}

\usepackage{amssymb}
\usepackage{amsmath,amssymb,amsfonts}
\usepackage{tikz}
\usepackage{tabularx}
\usepackage{natbib}
\usetikzlibrary{decorations.markings,positioning,shapes,arrows}
\usepackage{matlab-prettifier}

\newproof{proof}{Proof}

\begin{document}

\begin{frontmatter}

\title{An Efficient Numerical Method for an Approximate Solution of the Beam Equation}

\author[1]{Onur Baysal\corref{cor3}}
\ead{onur.baysal@um.edu.mt}
\author[2]{Maria Aquilina%
}
\ead{maria.aquilina.18@um.edu.mt}

\cortext[cor3]{Corresponding author}

\address[1]{Department of Mathematics, University of Malta, Malta}
\address[2]{Department of Mathematics, University of Malta, Malta}

\begin{abstract}
In this paper, we propose a horizontal type method of lines numerical scheme for the unsteady  Euler-Bernoulli beam equation. The problem is initially reformulated as a first order system of initial value problems and a suitable one-step difference scheme is used for the highest order temporal derivative which leads to a system of steady beam equations. Then resulted family of steady problems is solved iteratively by the finite element method with Hermite cubic basis functions.  This iterative procedure leads to approximations for both the solution of the unsteady problem and its derivatives. All these approximations are compared with the exact ones  to illustrate the performance of the proposed method. Moreover, the optimization of the mesh parameters is discussed  for both steady and unsteady problems by logarithmic scale plot. 
\end{abstract}

\begin{keyword}
Euler-Bernoulli Beam Equation, Method of Lines, Finite Element Method, Order of Accuracy.
\end{keyword}

\end{frontmatter}


\section{Introduction}
Solving a differential equation implies determining a certain quantity of a physical phenomenon. This can correspond to quantities such as temperature, pressure, velocity, displacement or density etc. This shows the importance of solving differential equations modelled by the relation of partial or ordinary derivatives in applied sciences. 

The finding of an analytical solution of a differential equation is often challenging, depending on the setting of the problem, and explicit solution may only be obtained for limited cases.
Numerical algorithms are considered a viable alternative to address these challenges. These algorithms can be classified into different categories based on various theoretical principles. Two of the most commonly used examples are the Finite Difference (FD) and Finite Element (FE) methods. Moreover, the combined use of these methods can result in highly effective hybrid techniques; for instance, the Method of Lines (MOL) approach may incorporate an appropriate combination of both. In this study, we propose a new hybrid numerical method for solving the unsteady Euler-Bernoulli beam equation.

Whether rigidly clamped, cantilevered, or simply supported, beams are integral components of many classical engineering systems, including airplane wings, stabilizers, building structures, bridge models, and helicopter rotor blades (see \cite{VIN2007} and references therein).
Recently, the source and boundary identification problems related to beam models have  emerged in fields such as medical diagnostics and nanoscale measurement systems, including transverse dynamic force microscopes (TDFM) \cite{ANTOG:2000}, \cite{NGU15} and atomic force microscopes (AFM) \cite{GB:XX2016}.

Despite its widespread use in applications, the demand for numerical solutions to this equation has significantly increased, yet there are only a few alternatives available in the literature. For a stationary model problem, the Finite Element (FE) approach is straightforward, utilizing cubic Hermite basis functions. \cite{WK&HB00}, (\cite{SG&DC&KJ&AS}). In the case of time-dependent problems, two examples of hybrid schemes given in \cite{WK&HB00} have been proposed based on the vertical Method of Lines (VMOL) approach. Essentially, VMOL is founded on the principle of independent discretization of spatial and temporal variables. More specifically, a semi-analytical structure is obtained by expressing the variational formulation within a finite-dimensional space, after which fully discrete algebraic system of equations are derived by substituting temporal derivatives with a proper finite difference expressions.
A similar but improved version, which offers enhanced stability in temporal discretization, has been applied to a different form of the beam equation, demonstrating both the efficiency and accuracy of the method \cite{AH&AK&OB24}.
Moreover, there have been numerous successful applications of this approach in prior studies, particularly for inverse problems \cite{OB&AH&AK24},\cite{AH&OB16}, and \cite{AH&OB&HI:19}. 
In this paper, we aim to develop a new method based on the horizontal Method of Lines (HMOL), which is not only easier to adapt to complex problem settings but also facilitates stable computation of approximate derivatives compared to the aforementioned methods. Unlike VMOL, HMOL begins with discretization in time, followed by an iterative solution of the resulting family of steady problems using the Finite Element Method (FEM). Furthermore, \cite{KSSJNR17} can be considered a comprehensive reference for space-time FE methods.

\begin{figure}[h]
  \centering
 \begin{tikzpicture}[scale=.8]

    \node[isosceles triangle,
	draw, rotate=90,
	fill=black!40,
	minimum size =.7cm] (T)at (0,-0.77){};
     \node[isosceles triangle,
	draw, rotate=90,
	fill=black!40,
	minimum size =.7cm] (T)at (4*3.14,-0.77){};

      \draw[dashed] (0,0) -- (4*3.14,0);
      \draw[thick,->] (4*3.14,0) -- (13.5,0) node [right]{$x$};
      \draw[thick,->] (0,-1.5) -- (0,1.8) node [left]{$w$};
      \draw[ultra thick,color=black] (4*3.14,-0.1) -- (4*3.14,0.1);
      \draw[color=black,ultra thick,smooth,domain=0:{4*3.14}] plot (\x,{0.1+sin(deg(\x)/2)/2});
      \draw[color=black,ultra thick,smooth,domain=0:{4*3.14}] plot (\x,{-0.1+sin(deg(\x)/2)/2});
      \node[label=right:{$w(0,t)=0$}] at (0,1.3) {};
      \node[label=right:{$w_{xx}(0,t)=0$}] at (0,0.7) {};
      \node[label=center:{$w(l,t)=0$}] at (4*3.14,1.3) {};
      \node[label=center:{$w_{xx}(l,t)=0$}] at (4*3.14,0.7) {};
  \end{tikzpicture}
\caption{{\small Geometry of the problem: Hinged-hinged beam model.}}
  \label{geometry}
\end{figure}
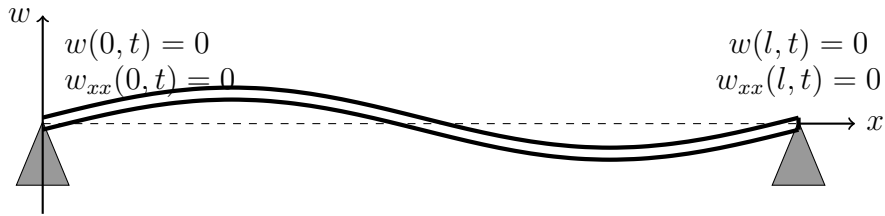

The non-homogenous dynamic Euler-Bernoulli beam with supported (hinged) boundary conditions is  illustrated in Figure \ref{geometry}, while the non-trivial boundary condition case  is modelled by the following equation:
\begin{eqnarray}\label{beam01}
\left\{\begin{array}{ll}
	\rho(x)w_{tt}+\eta(x)w_t+\mathcal{L}w=g(x,t)~~\text{for} ~~(x,t)\in(0,l)\times(0,T), \\
    w(x,0)=p(x), ~ w_t(x,0)=q(x) ~~\text{for} ~~ x\in (0,l),\\
	w(0,t)=a(t),~r(0) w_{xx}(0,t)=\widetilde{a}(t),~w(l,t)=b(t),~r(l)w_{xx}(l,t)=\widetilde{b}(t)~~\text{for} ~~ t\in[0,T].
\end{array}\right.
\end{eqnarray}

\noindent Here $\mathcal{L}w:=(r(x)w_{xx})_{xx}-(\mathfrak{j}(x)w_{x})_x+s(x)w$ and  $w(x,t)$  is the displacement function, depending on the space $x\in(0,l)$  and time $t \in[0,T]$  variables with $l,T\in \mathbb{R}^+$. In addition,  $g(x,t)$ is the load distribution, $r(x)=EI(x)$ is the flexural rigidity where $E>0$  is the elasticity modulus and  $I(x)>0$ is the moment of inertia of the cross-section, $\rho(x)$ is the mass density of the beam, $\eta(x)$ and $\mathfrak{j}(x)$ denote the damping coefficient and the traction force, respectively.\\
We also assume the following conditions:
\begin{eqnarray} \label{cond1}
\left \{ \begin{array}{ll}
\rho, r, \eta, \mathfrak{j} \in L^{\infty}(0,l)~~ \text{and}~~p \in H^{0,2}(0,l), ~ q \in L^2(0,l)\\
0<\underline{\rho}\leq\rho(x)\leq \overline{\rho}~~ \text{and}~~ 0< \underline{r} \leq r(x) \leq \overline{r}~~\text{for all}~~ x\in(0,l).
\end{array} \right.
\end{eqnarray}
Here the notation $H^k(0,l)$ is used for the Hilbert-Sobolev space with order $k$ and $H^{0,k}(0,l):=H^k(0,l)\cap C_0(0,l).$

This paper is organised as follows. In Section 2, the steady linear beam equation and corresponding weak formulation are  introduced. In addition, a simple FE method with cubic basis functions is developed and its performance analysis illustrated with a simple test problem.  Section 3 is dedicated to unsteady problem were the main algorithm  based on the HMOL approach constructed. All details of this method are given and its performance is analyzed using a test problem with different mesh parameters. Through this analysis, the optimization of the mesh parameters is discussed. 

\section{Stationary Problem}
In this section, an efficient  and simple numerical algorithm is considered for the following boundary value problem related to stationary beam model.
\begin{eqnarray}\label{beam_steady_n}
\left\{\begin{array}{ll}
	 (r(x)w'')''+s(x)w-(\mathfrak{j}(x)w')'=f(x),~~~ \text{for}~~x\in (0,l), \\
	w(0)=a, ~~ r(0)w''(0)=\widetilde{a}, ~w(l)=b, ~~r(l)w''(l)=\widetilde{b}.
\end{array}\right.
\end{eqnarray}

We assume that the inputs in (\ref{beam_steady_n}) satisfy the following basic conditions:
\begin{eqnarray} \label{beam_steady_c}
\left \{ \begin{array}{ll}
 f\in L^2(0,l)~~ \text{and} ~~ s,\mathfrak{j} \in L^\infty(0,\ell)\\ [3pt]
 r \in L^\infty(0,l)~~ \text{with}~~  0< \underline{r} \leq r(x) \leq \bar{r}
\end{array} \right.
\end{eqnarray}
Then weak formulation of (\ref{beam_steady_n}) is as follows. Find $w\in \mathbb{V}$ such that, for all $ v\in \mathbb{V}_0$ following holds,
\begin{eqnarray}\label{weak_steady}
(r(\cdot)w'',v'')+(s(\cdot)w,v)+(\mathfrak{j}(\cdot)w',v')=(f,v)+\widetilde{b}v'(l)-\widetilde{a}v'(0).
\end{eqnarray}
Here $\mathbb{V}:=\{w\in H^2(0,l): w(0)=a, ~~w(l)=b\}$ and $\mathbb{V}_0:=\{w\in H^2(0,l): w(0)=w(l)=0\}.$\\

\noindent \textbf{Remark:} Assuming homogenous boundary condition is not a restriction. One can easily obtain the problem in (\ref{beam_steady_n}) with homogenous boundary conditions by using following cubic auxiliary function which satisfies $\theta(0)=a,\ r(0)\theta''(0)=\widetilde{a},\ \theta(l)=b,\ r(l)\theta''(l)=\widetilde{b}$.
\begin{eqnarray*}\label{zxt}
       \theta(x) = a\dfrac{(l-x)}{l} + b \dfrac{x}{l} + \frac{\widetilde{a}}{r(0)} \biggl[\dfrac{(l-x)^{3}}{6l} - \dfrac{l}{6} (l-x)\biggr] + \frac{\widetilde{b}}{r(l)} \biggl[\dfrac{x^{3}}{6l} - \dfrac{xl}{6}\biggr]
\end{eqnarray*}
Then the solution of (\ref{beam_steady_n}) can be obtained by taking $w=\widetilde{w}+\theta$ where
\begin{eqnarray}\label{beam_steady}
\left\{\begin{array}{ll}
	 (r(x)\widetilde{w}'')''+s(x)\widetilde{w}-(\mathfrak{j}(x)\widetilde{w}')'=\widetilde{f}(x),~~~ \text{for}~~ x\in (0,l), \\
	\widetilde{w}(0)= r(0)\widetilde{w}''(0)=\widetilde{w}(l,t)=r(l)\widetilde{w}''(l)=0.
\end{array}\right.
\end{eqnarray}
with $\widetilde{f}(x)=f(x)-(r(x)\theta(x))''-s(x)\theta(x)+(\mathfrak{j}(x)\theta'(x))'$. This trick provides some advantages in theoretical analysis.
One can prove that the weak problem defined in (\ref{weak_steady}) with $a=b=\widetilde{a}=\widetilde{b}=0$ has a unique solution by the Lax-Milgram theorem (\cite{BH11}) by using the following $\mathbb{V}$-norm in $H^{2}(0,l)$.
  $$||w||_\mathbb{V}=\sqrt{||w_{xx}||^2+||w||^2} $$

\subsection{Numerical Solution of the Steady Problem by the Finite Element Method}\label{sec2.1}

This section will review some fundamental ideas for the unsteady state. We do not claim to propose any new approach for this steady case but it contains a detailed receipt of computer programming and the determination of an optimum mesh parameter. We begin with the Galerkin approximation of the problem (\ref{weak_steady}). Let $0=x_1< x_2< \cdots <x_{M+1}=l$ be a partition of $(0,l)$ with the element size $|e_m|=h_{m}=x_{m+1}-x_{m}$ for $m=1:M+1$. Assume
$\mathbb{V}_h=span\ \{\phi_m\}_{m=1}^{2(M+1)}$  is a finite dimensional subspace of $H^{2}(0,l)$ whose elements are Hermite cubic basis functions such that
$$\phi_n|_{e_m}(x)=\left \{ \begin{array}{ll}
1-3\bar{x}^2/h_m^2+2\bar{x}^3/h_m^3  & \text{if} ~~ n=2\ m-1 \\ [3pt]
\bar{x}-2\ \bar{x}^2/h+\bar{x}^3/h^2 & \text{if} ~~ n=2\ m\\[3pt]
3\bar{x}^2/h_m^2-2\bar{x}^3/h_m^3  & \text{if} ~~ n=2\ m+1\\ [3pt]
-\bar{x}^2/h+\bar{x}^3/h^2 & \text{if} ~~ n=2\ (m+1)\\ [3pt]
0 & \text{otherwise}.
\end{array} \right.
$$
where $\bar{x}=x-x_m$. Then the Galerkin Finite Element Method related to (\ref{weak_steady}) can be expressed as follows.\\

\noindent  Find $w_h=\sum_{m=1}^{2(M+1)}\alpha_m\phi_m\in \mathbb{V}_h$ such that, for all $ v_h\in \mathbb{V}^0_h:=\mathbb{V}_h\cap H^{0,2}(0,l)$ the following holds,
\begin{eqnarray}\label{Fem_steady}
\left \{ \begin{array}{ll}(r(\cdot)w_h'',v_h'')+(s(\cdot)w_h,v_h)+(\mathfrak{j}(\cdot)w_h',v_h')=(f,v_h)+\widetilde{b}v_h'(l)-\widetilde{a}v_h'(0).\\
w_h(0)=\alpha_1=a, ~~ w_h(l)=\alpha_{2M+1}=b,
\end{array} \right.
\end{eqnarray}

By replacing test functions $v_h$ with basis function $\phi_m$, the following $(2M+2)\times(2M+2)$ system of algebraic equations is obtained.

\noindent Find $\Gamma=(\alpha_1,\ \alpha_2, \cdots \alpha_{2(M+1)})$ such that,
$R\ \Gamma = C.$ \\

\noindent Here the components of the $(2M+1)$ by $(2M+1)$ matrix $R$ and the right hand side vector $C$ are as follows:
\begin{eqnarray}\label{Components} \begin{array}{ll}
R_{ij}=(r(\cdot)\phi_j'',\phi_i'')+(\mathfrak{j}(\cdot)\phi_j',\phi_i')+(s(\cdot)\phi_j,\phi_i)~~ \text{for}~~{i,j}=1:2M+2, \\
C_{j}=(f,\phi_j)+\widetilde{b}\phi_j'(l)-\widetilde{a}\phi_j'(0)~~ ~~ \text{for}~~{j}=1:2M+2.
\end{array}
\end{eqnarray}
The calculation of these components in (\ref{Components}) is performed at the element level, which is a routine procedure for executing any FEM algorithm. Moreover, these integrals can be efficiently evaluated using Gauss quadrature formulas \cite{ES&DM12}. An effective application of this type of calculation can be found in \cite{AM12} for reduced HCT-basis functions related to the plate equation. In this study, the same method is adapted for use with the beam equation. An example of MATLAB code for element calculations related to the constant elasticity term $r\in \mathbb{R^+}$, with a step size $h$, is shown below.\\

\begin{lstlisting}[
frame=single,
%numbers=left,
style=Matlab-editor]
 pt=[-sqrt(3/5) 0 sqrt(3/5)]; wg=[5/9 8/9 5/9];
 for i=1:3
   D2=[6*pt(i)/h; 3*pt(i)-1; -6*pt(i)/h; 3*pt(i)+1]/h;
   R_el=R_el+h/2*wg(i)*r*(D2*D2');
 end
\end{lstlisting}

\noindent Note that, due to the given essential boundary conditions $\alpha_1=a$ and $\alpha_{2M+1}=b$,
we simplify the system $R\ \Gamma = C$ by using the following corrections on $R$ and $C$.
\begin{eqnarray}\label{Corection} \begin{array}{ll}
R_{1j}=\delta_{1,j}~~\text{and}~~ R_{2M+1j}=\delta_{2M+1,j} ~~ \text{for}~~{j}=1:2M+2, \\
C_{1}=a ~~\text{and}~~ C_{2M+1}=b.
\end{array}
\end{eqnarray}

\subsection{Test Problem for the Steady Beam}

In this section, we demonstrate the performance and accuracy of the proposed algorithm for the weak formulation (\ref{Fem_steady}) using the following parameters:
 $$r(x)=1+x,~~ s(x)=\cos(x),~~\mathfrak{j}(x)=3$$
 \noindent with domain and discretization parameters of $l=1$ and $h=1/50$, respectively, the right-hand side function $f(x)$ and boundary conditions can be derived for $w(x)=1-x+\sin^2x$ to serve as an exact solution.

\begin{figure}[!htb] \centering	\includegraphics[width=6.75cm,height=6.25cm]{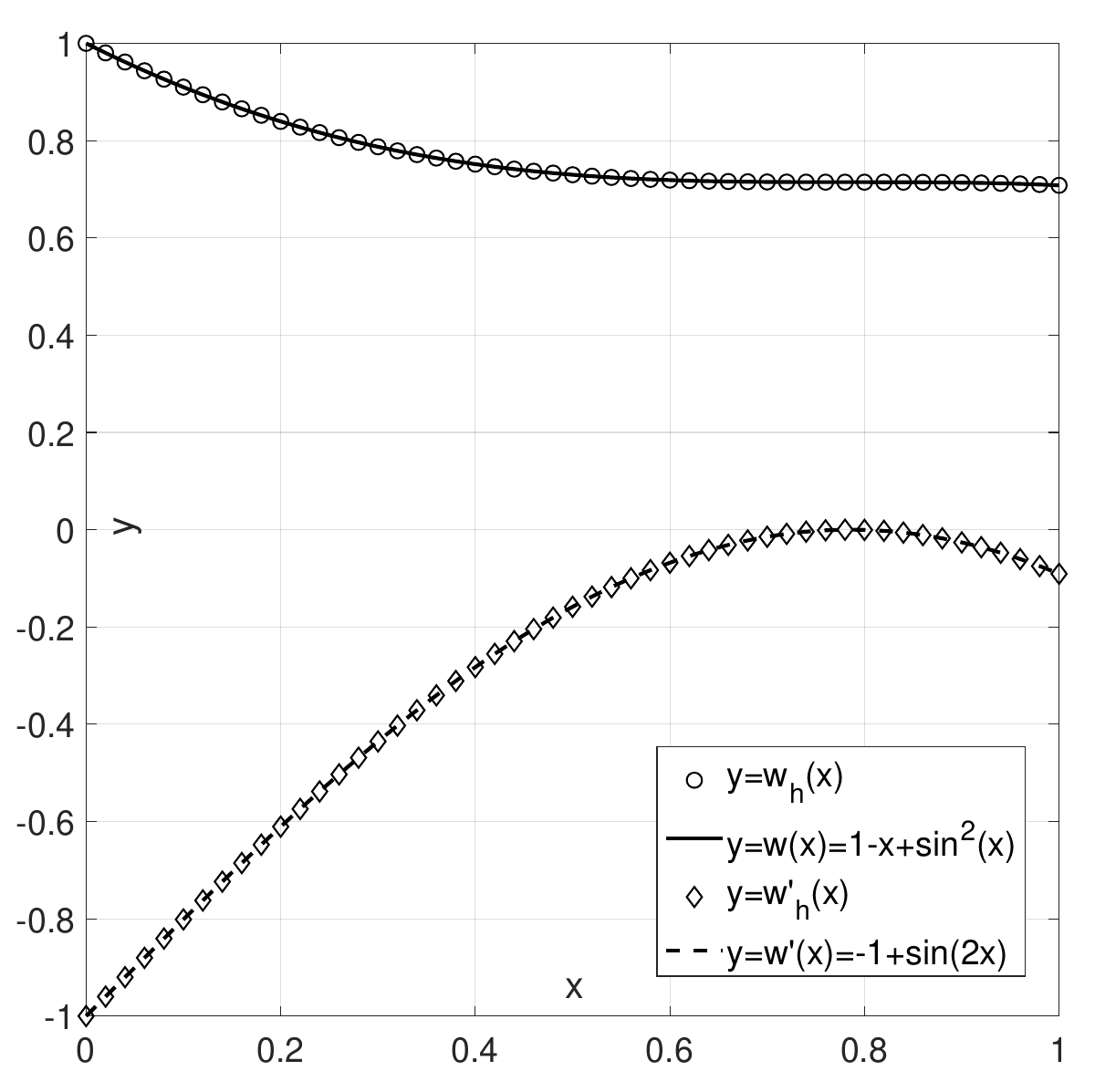}~~\includegraphics[width=6.75cm,height=6.25cm]{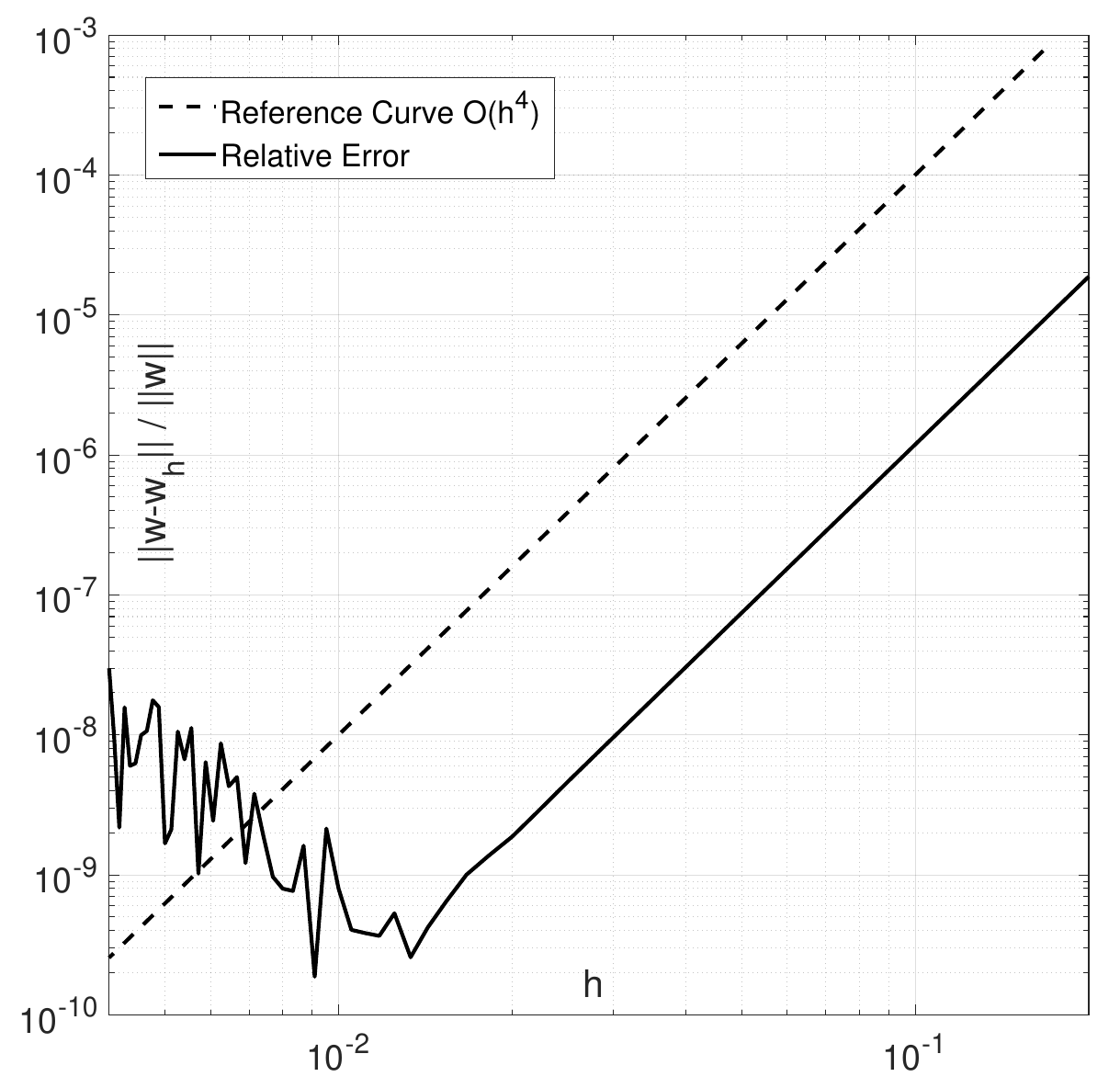}
	\caption{Comparison of the exact solution and approximate solution with derivatives (left), error plot in log scale  (right).} \label{fig_steady}
\end{figure}
Figure 2 illustrates the accuracy analysis of the method. The plot on the right-hand side shows that the accuracy error is of order 4, provided that the step size range is approximately $h \in [1/80, 1/4]$. Additionally, to simplify the notation, we omit the subscript $L^2(0,l)$ in $||\cdot||_{L^2(0,l)}$  in all figures.

\section{Numerical Solution of the Dynamical Beam Equation} \label{NSDP}

In this section, we introduce a new numerical method for solving the damped unsteady beam equation defined in (\ref{beam01}), assuming $\rho\equiv 1$. The existence and uniqueness of the weak solution to (\ref{beam01}) can be established using the Galerkin method. The proof for the clamped case, as discussed in \cite{OB&AH19} and \cite{AH&AGR2021}, can be easily adapted to the hinged-hinged case, following the principles of this technique detailed in \cite{ELC02}.

To simplify presentation, we first define $z=w_t$, and rewrite problem (\ref{beam01}) as follows.
\begin{eqnarray}\label{vector}
\left(
    \begin{array}{c}
      w_t \\
      z_t \\
    \end{array}
  \right)=
\left(
    \begin{array}{c}
      z \\
      -\eta z - \mathfrak{L}w +g \\
    \end{array}
  \right)~~ \text{with} ~~
\left(
    \begin{array}{c}
      w(\cdot,0) \\
      w_t(\cdot,0) \\
    \end{array}
  \right)=
\left(
    \begin{array}{c}
      p \\
      q \\
    \end{array}
  \right)
 \end{eqnarray}
 This is a routine expression for partial differential equations (PDEs) that involve higher-order temporal derivatives, such as the wave equation \cite{MA20}. The first-order initial value problem represented by the two equations given in (\ref{vector}) can be discretized using a suitable explicit or implicit time integrator on the uniform partition $\{t_n\}_{n=1}^{N+1}$ of $[0,T]$ with a step size $\tau:=T/N$. For this operation, we utilize the Crank-Nicolson (trapezium) method which is easily adapted, A-stable, and second-order accurate. The resultant temporal discrete form related to (\ref{vector}) is as follows
 
\begin{eqnarray*}\label{discvector}
\left(
    \begin{array}{c}
      W_{n+1} \\
      Z_{n+1} \\
    \end{array}
  \right)= \left(
    \begin{array}{c}
      W_{n} \\
      Z_{n} \\
    \end{array}
  \right)
+\displaystyle \frac{\tau}{2}
\left(
    \begin{array}{c}
      Z_n+Z_{n+1} \\
      -\eta (Z_n+Z_{n+1}) - \mathfrak{L}(W_n+W_{n+1}) +G_n+G_{n+1} \\
    \end{array}
  \right)
 \end{eqnarray*}
\noindent with  $W_1=p$ and $Z_1=q$ (given initial conditions). Here $G_n:=g(\cdot,t_n),\ W_{n}\approx w(\cdot,t_n)$ and $Z_{n}\approx w_t(\cdot,t_n)$ for $n=1:N+1$. Then, we obtain the following couple of equations.
\begin{eqnarray}\label{discvector1}
\left\{\begin{array}{ll}\displaystyle
	 Z_{n+1}=\frac{2}{\tau}(W_{n+1}-W_n), \\ [5 pt]
	(r(x)\ W_{n+1}'')''+\widetilde{s}(x)\ W_{n+1}-(\mathfrak{j}(x)\ W_{n+1}')'=\widetilde{g}(x)
\end{array}\right.
\end{eqnarray}
with the boundary conditions 
\begin{eqnarray}\label{discBC}
W_n(0)=a_n,~~ r(0)W_n''(0)=\widetilde{a}_n,~~ W_n(l)=b_n,~\text{and}~ r(l)W''_n(l)=\widetilde{b}_n 
\end{eqnarray}
 where $a_n:=a(t_n),~\widetilde{a}_n:=\widetilde{a}(t_n),~ b_n:=b(t_n),~\widetilde{b}_n:=\widetilde{b}(t_n) $ and
\begin{eqnarray}\label{discvector2}
\left\{\begin{array}{ll}\displaystyle
	\widetilde{s}(x)=\frac{4}{\tau^2}+\frac{2}{\tau}\ \eta(x)+s(x), \\
	\displaystyle \widetilde{g}(x)=-(r(x)\ W_{n}'')''+\left(\frac{4}{\tau^2}+\frac{2}{\tau}\eta(x)+s(x)\right)W_n+\frac{4}{\tau}Z_n+(\mathfrak{j}(x)\ W_n')'+G_n+G_{n+1}.
\end{array}\right.
\end{eqnarray}
The iteration prescribed  in (\ref{discvector1}) with (\ref{discvector2})  is well-defined. From the given $W_1$ and $Z_1$, we define their projections $W_{h,1}$ and $Z_{h,1}$ on $\mathbb{V}_h$ such that
$$(W_1-W_{h,1},v_h)=0 ~~ \text{and} ~~ (Z_1-Z_{h,1},v_h)=0~~  \text{for all} ~~v_h\in \mathbb{V}_h.$$
Then one can compute  $W_{h,2}\approx W_{2}$ using the second part of (\ref{discvector1}) and the numerical method provided for the steady equation in the previous section. Subsequently, $Z_{h,2}\approx Z_{2}$ can be obtained from the first part of (\ref{discvector1}). This process is continued for $N$ steps, resulting in the series of functions $W_{h,n}$ and $Z_{h,n}$ for $n=1:N+1$.

Note that, in the finite element calculation, $(\widetilde{g}(\cdot),v_h)$ can not be computed directly due to the term
$(r(\cdot)W_{h,n}'')''$. It is known that $W_{n}\in \mathbb{V}_h$ so $(r(\cdot)W_{h,n}'')'',v_h)$ can not be $L_2$ integrable on $(0,l)$. Instead, we apply integration by parts to obtain
$$((r(\cdot)W_{h,n}'')'',v_h)=(r(\cdot)W_{h,n}'',v_h'')-\widetilde{b}_nv_h'(l)+\widetilde{a}_nv_h'(0).$$
Similarly, the following identity simplifies the corresponding finite element calculation.
$$((\mathfrak{j}(\cdot)\ W_{h,n}')',v_h)= -(\mathfrak{j}(\cdot)\ W_{h,n}',v_h')$$

\subsection{Test Problem for the Dynamical Beam Equation}

Here, the performance analysis of the proposed method in section \ref{NSDP}  is presented by setting the following parameters in (\ref{beam01}):
$$r(x)=1+x, ~~ \eta(x)=2, ~~\mathfrak{j}(x)=3 x,~~ s(x)=\cos\ x.$$
Assuming the exact solution to be $w(x,t)=(1-x+\sin^2x)(\sin t)$ for $x\in[0,1]$ and $t\in[0,1]$, one can compute the right hand side function $g(x,t)$, the initial conditions  $p(x),~q(x)$ for $x\in[0,1]$  and the boundary conditions $a(t),~\widetilde{a}(t),~\ b(t),~\widetilde{b}(t)$ for  $t\in[0,1]$.\\
To compare the approximate and analytical solutions, we select constant spatial and temporal mesh parameters $h=1/50$ and $\tau=1/1250$, respectively, and define ${W_h(x,t)}$ as the linear interpolation of the set of all solutions  $\{W_{h,n}\in V_h\}_{n=1}^{N+1}$ in the temporal dimension such that  for $n=1,\, \cdots,\, N$,
$$W_h(x,t)|_{[t_n,t_{n+1}]}:=\frac{t-t_{n}}{\tau}W_{h,n+1}(x)-\frac{t-t_{n+1}}{\tau}W_{h,n}(x)\approx w(x,t)|_{[t_n,t_{n+1}]}.$$
In a similar way, $Z_h(x,t)\approx w_t(x,t)$ can be defined. Moreover, we denote the approximate slopes of deflection and velocity, namely $\partial W_h/\partial x $ and $\partial Z_h/\partial x $  with $\partial_x W_h$ and $\partial_x Z_h$, respectively. These functions are directly obtained from (\ref{Fem_steady}) and (\ref{discvector1})-(\ref{discvector2}) by utilizing  Hermite cubic basis functions.

\begin{figure}[!htb]	 \centering \includegraphics[width=7cm,height=6.5cm]{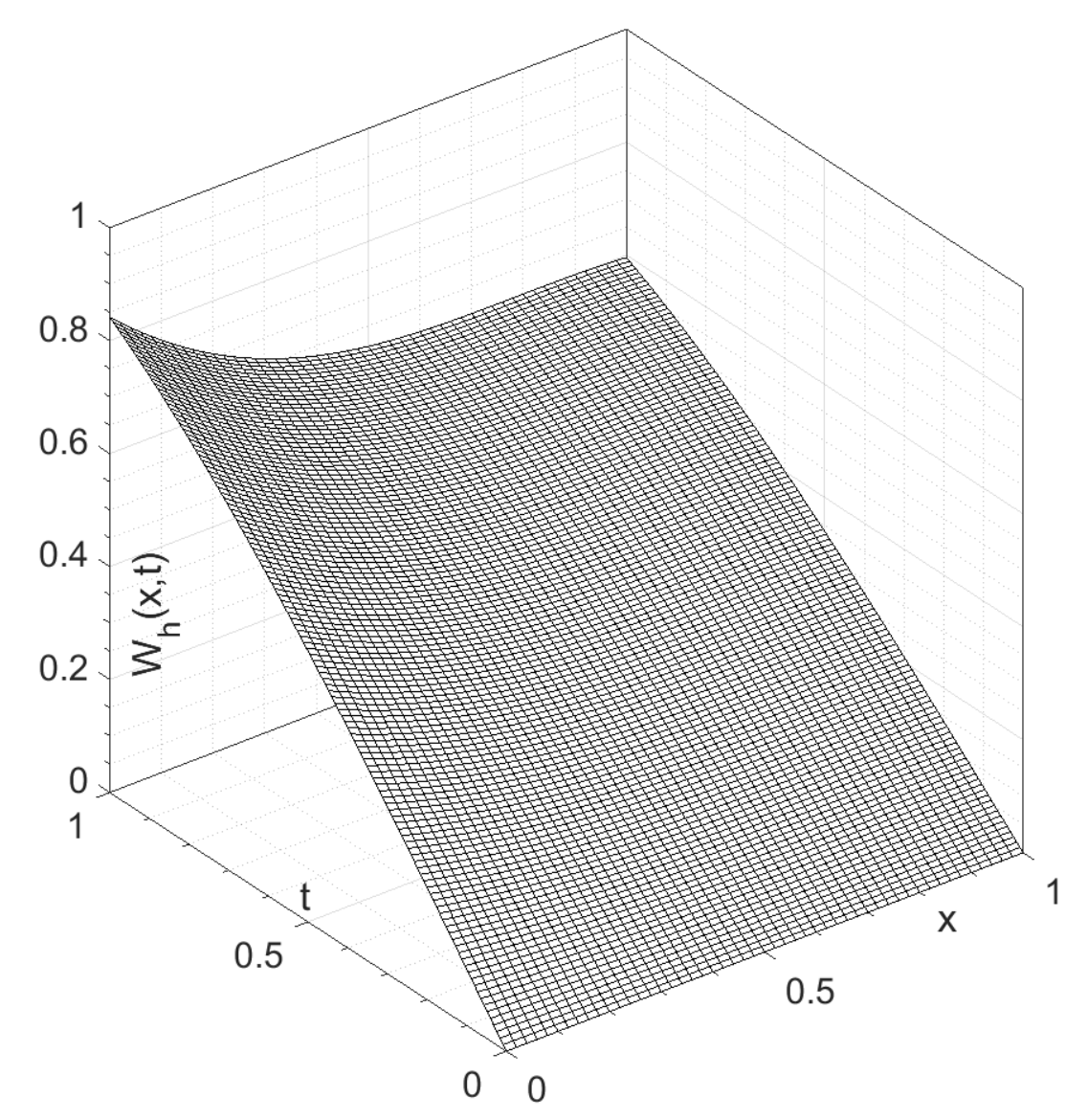}\includegraphics[width=7cm,height=6.5cm]{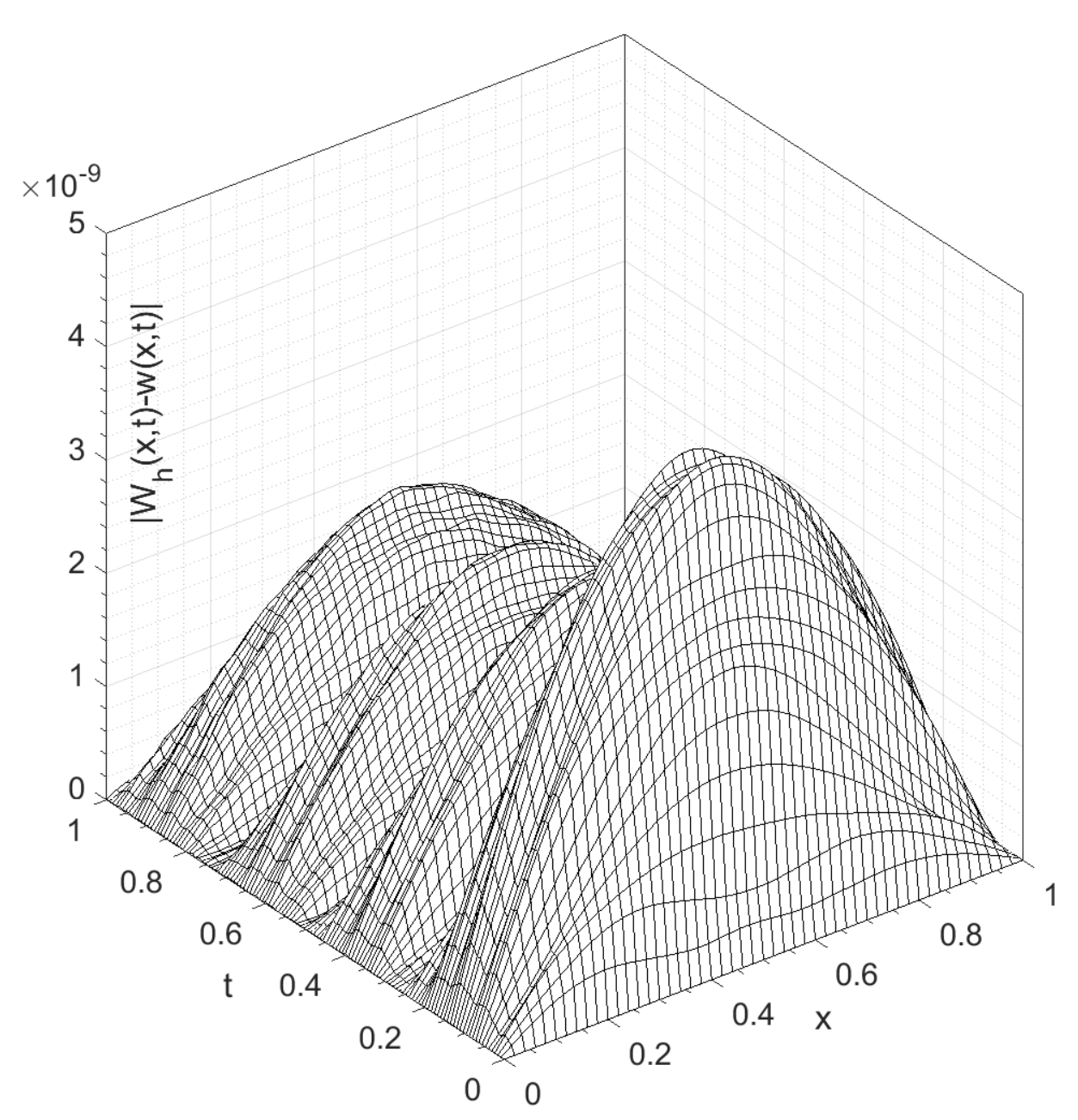}
	\caption{Approximate solution $W_h(x,t)$ (left) and corresponding error (right).}\label{fig_U}
\end{figure}

\begin{figure}[!htb]	 \centering \includegraphics[width=7cm,height=6.5cm]{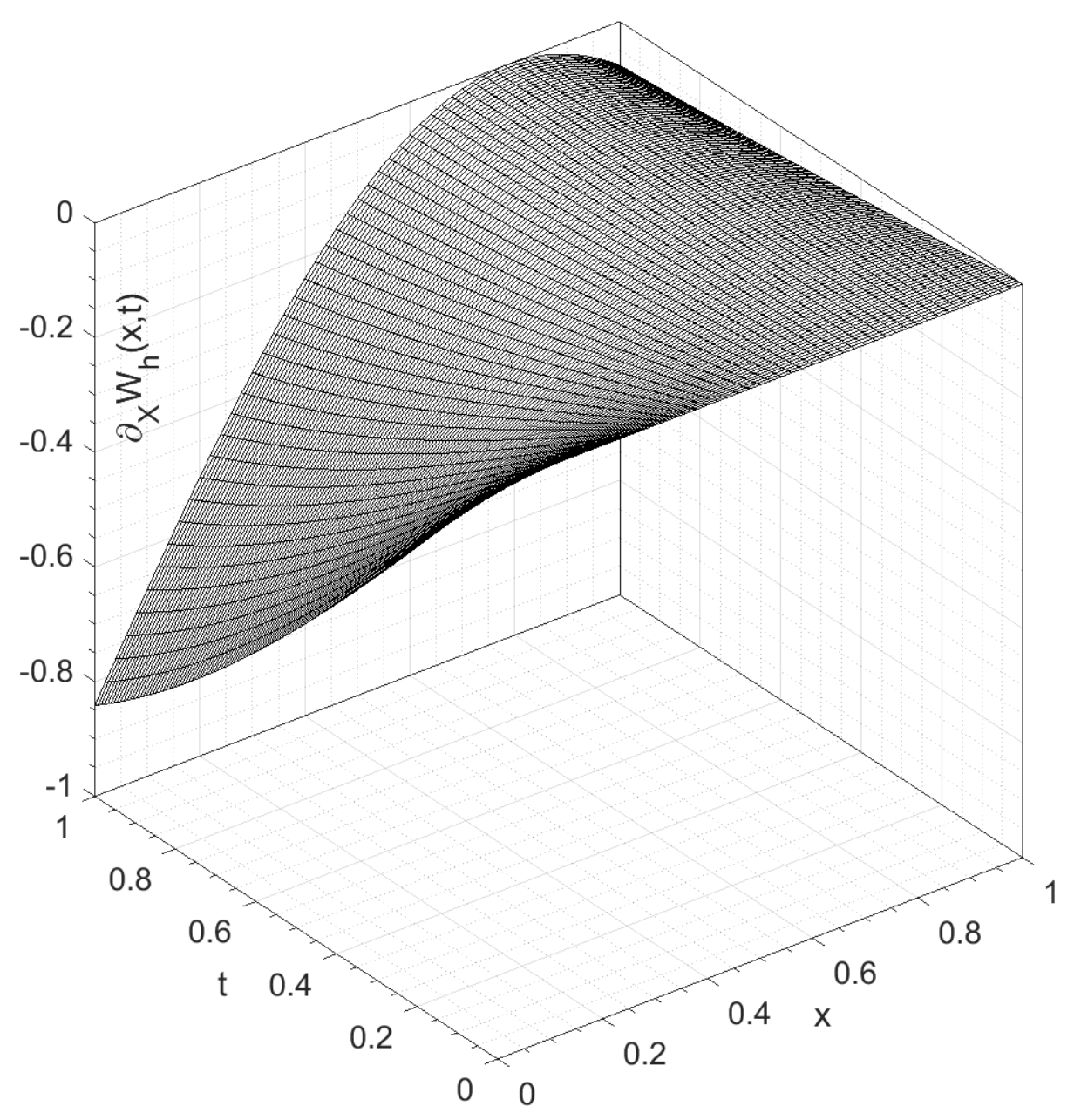}\includegraphics[width=7cm,height=6.5cm]{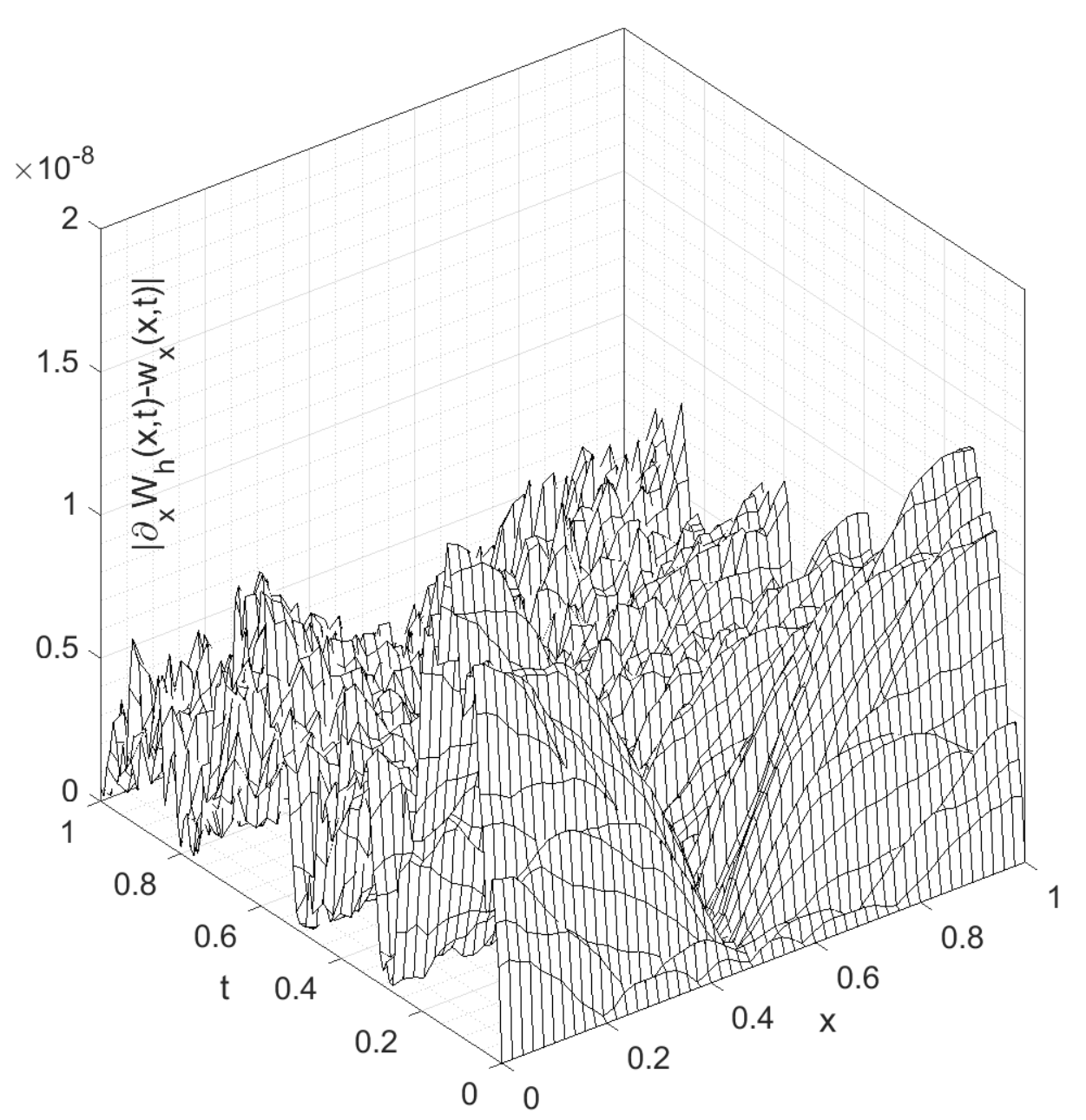}\label{fig_Ux}
	\caption{Approximate solution $\partial_xW_h(x,t)$ (left)  and corresponding error (right).}
\end{figure}

\begin{figure}[!htb]	 \centering \includegraphics[width=7cm,height=6.5cm]{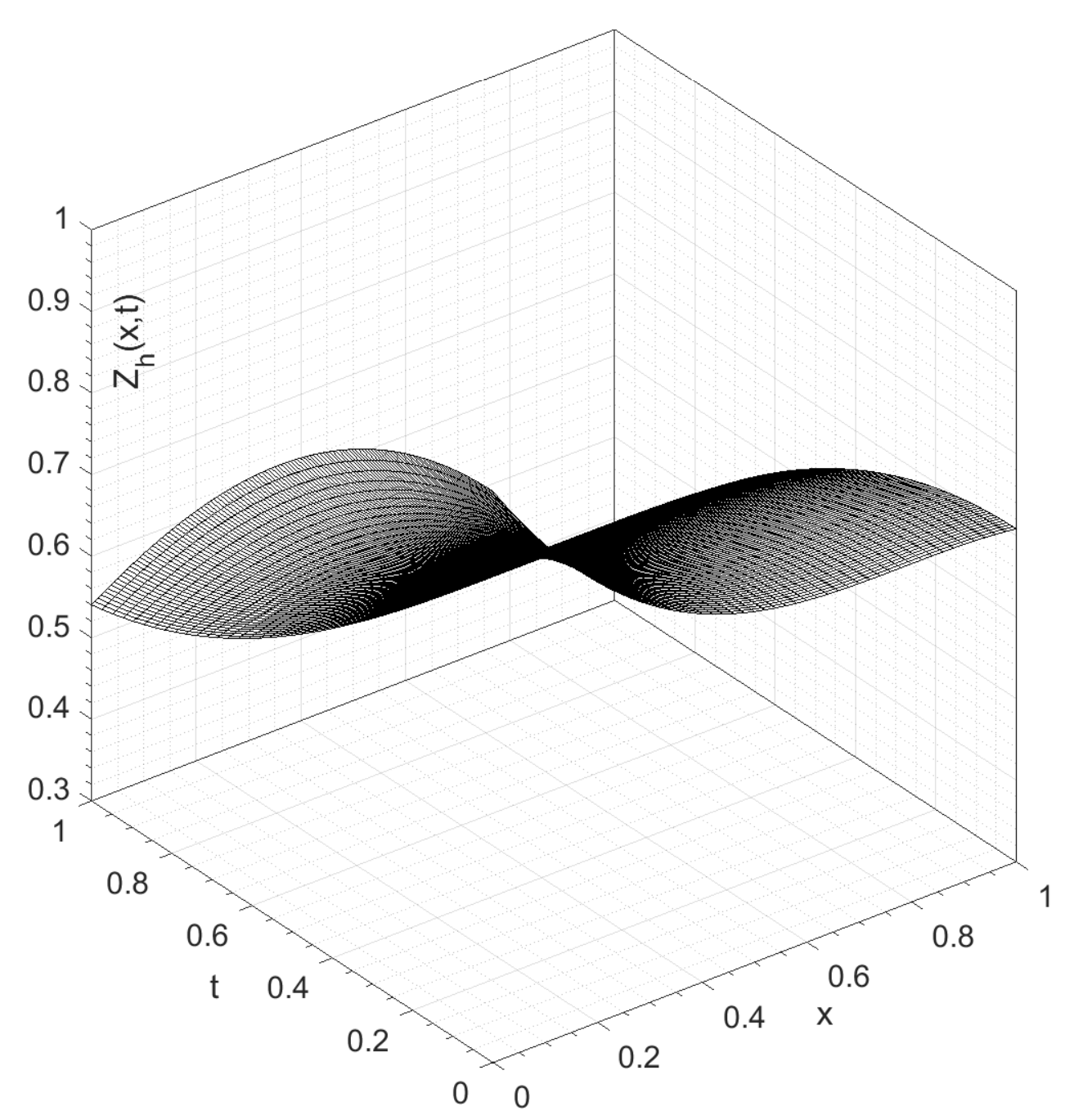}\includegraphics[width=7cm,height=6.5cm]{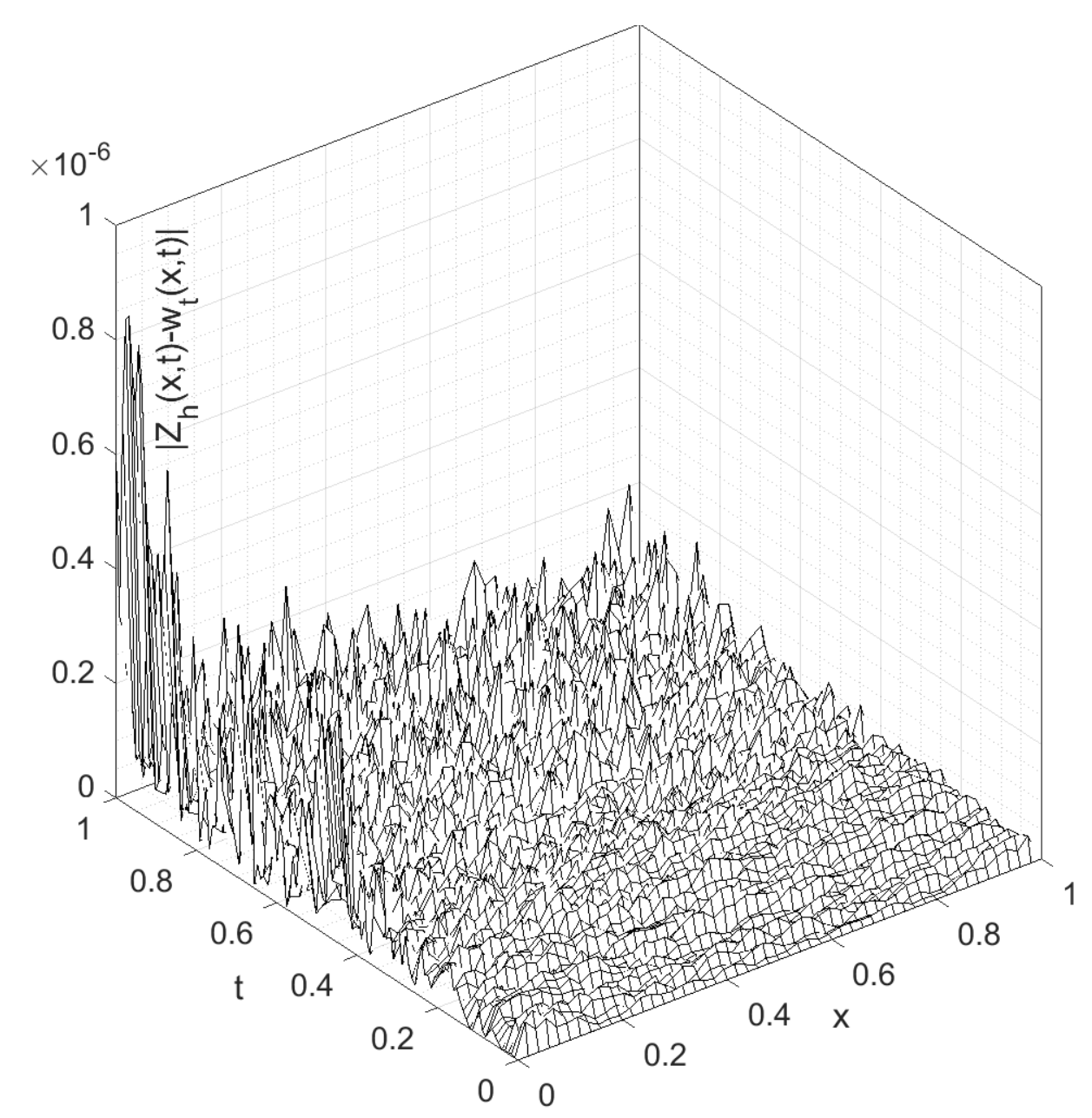}
	\caption{Approximate solution $Z_h(x,t)$ (left) and corresponding error (right).}\label{fig_V}
\end{figure}

\begin{figure}[!htb]	 \centering \includegraphics[width=7cm,height=6.5cm]{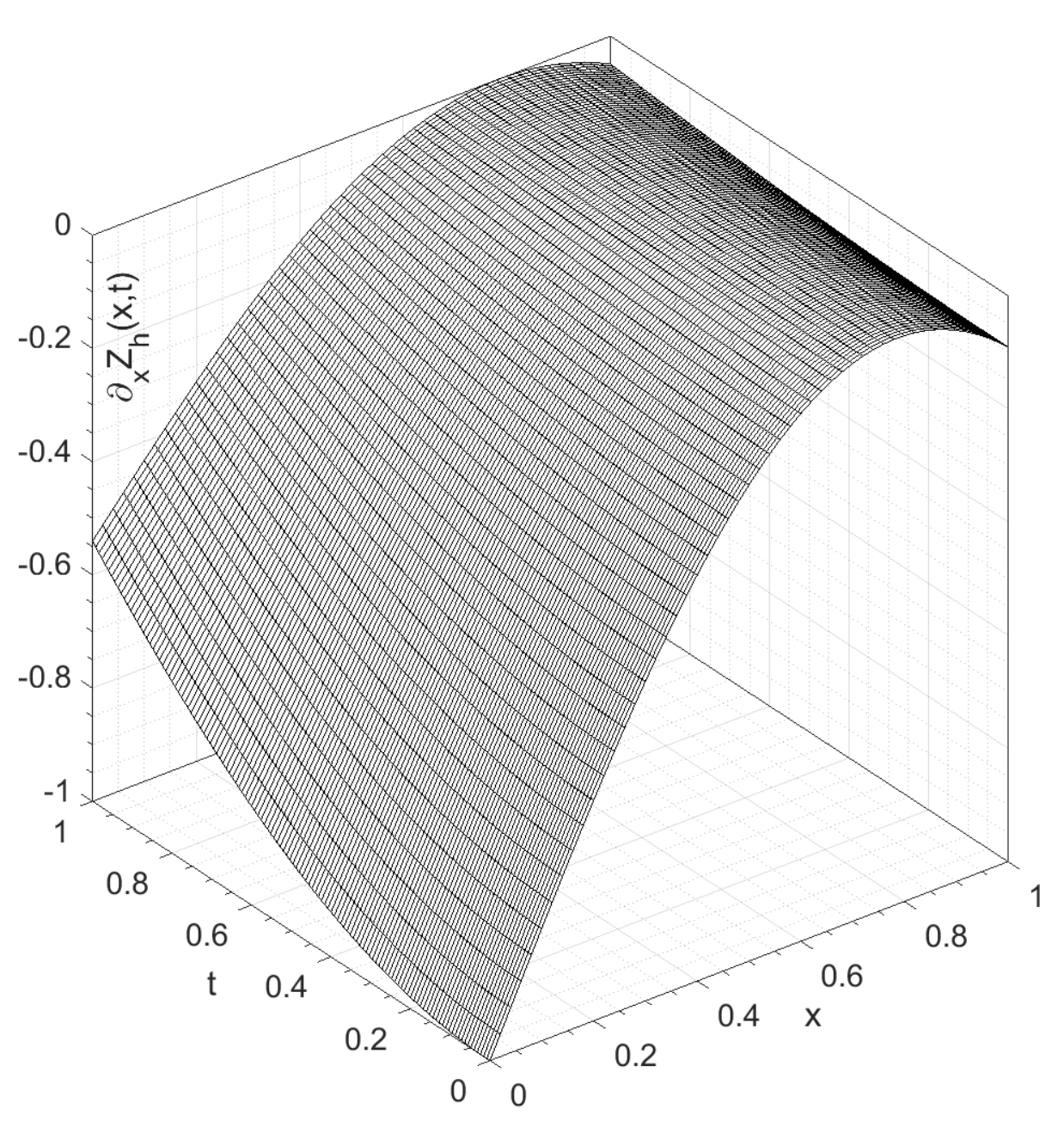}\includegraphics[width=7cm,height=6.5cm]{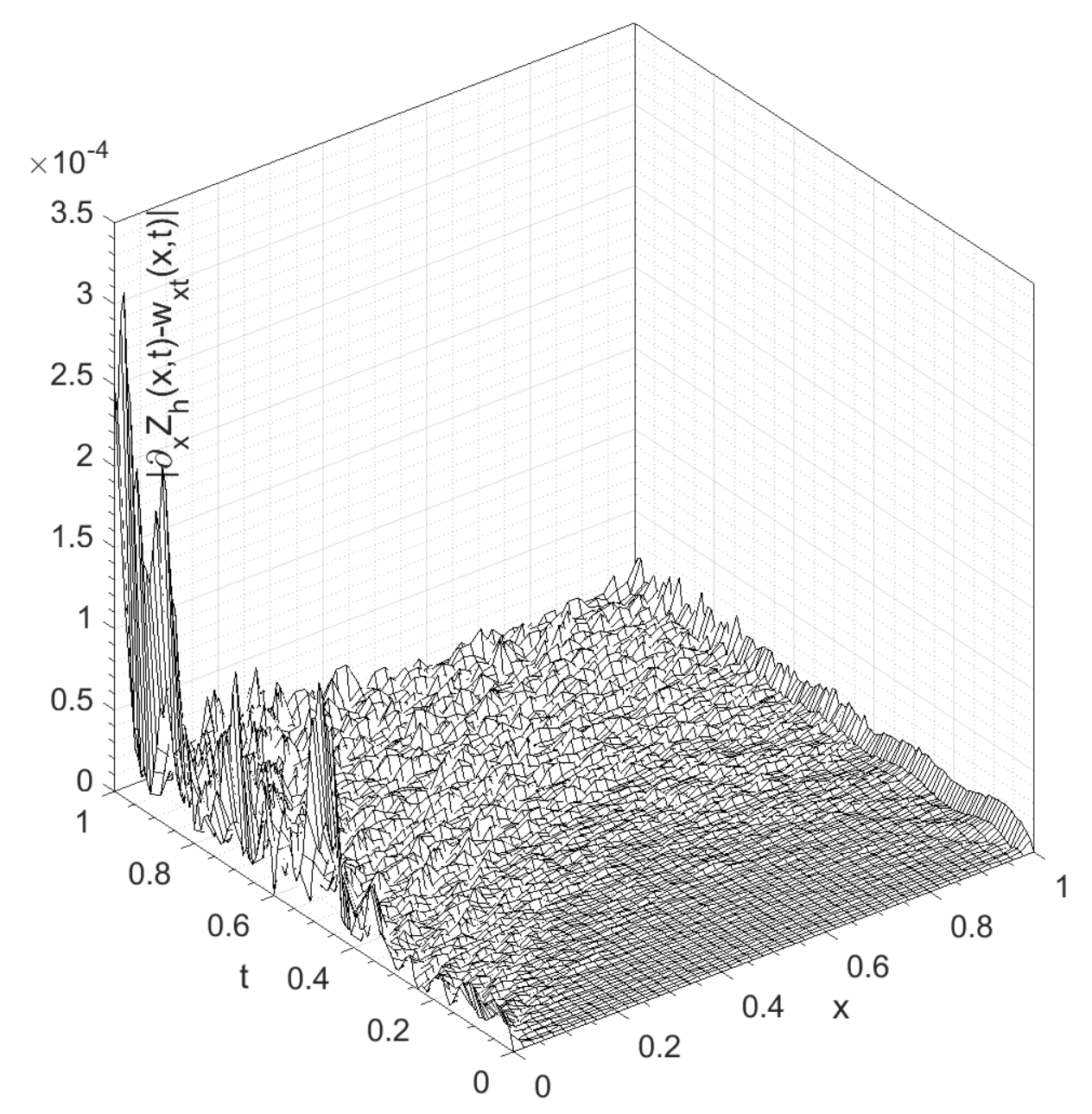}
	\caption{Approximate solution $\partial_xZ_h(x,t)$ (left)  and corresponding error (right).}\label{fig_Vx}
\end{figure}

The efficiency of the proposed iterative scheme based on the HMOL approach is illustrated in Figs. (\ref{fig_U}) and (\ref{fig_Ux}) for $W_h(x,t)\approx w(x,t)$ and $\partial_xW_h(x,t)\approx w_x(x,t)$, respectively. It is well known that the evaluation of numerical derivatives can be quite challenging due to computational errors. This difficulty increases when capturing higher order derivatives, often making it nearly impossible to obtain stable results.
However, the proposed algorithm  approximates accurately not only $w_t(x,t)$ but also $w_{xt}(x,t)$, as shown in Figures (\ref{fig_V}) and (\ref{fig_Vx}). This potential can be regarded as one of the most significant advantages of the algorithm.

\begin{figure}[!htb]	 \centering \includegraphics[width=6.75cm,height=6.28cm]{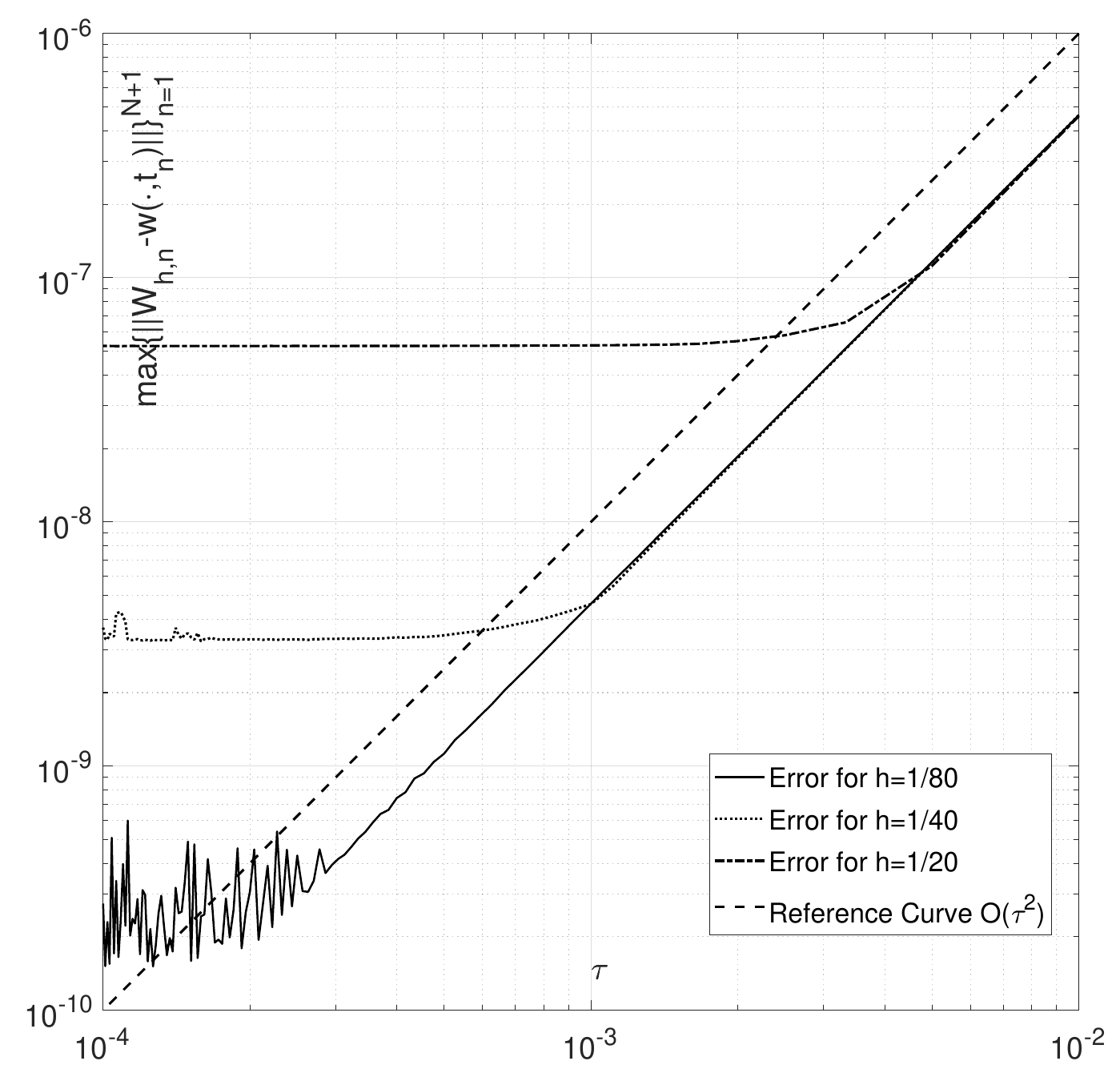}~~\includegraphics[width=6.75cm,height=6.25cm]{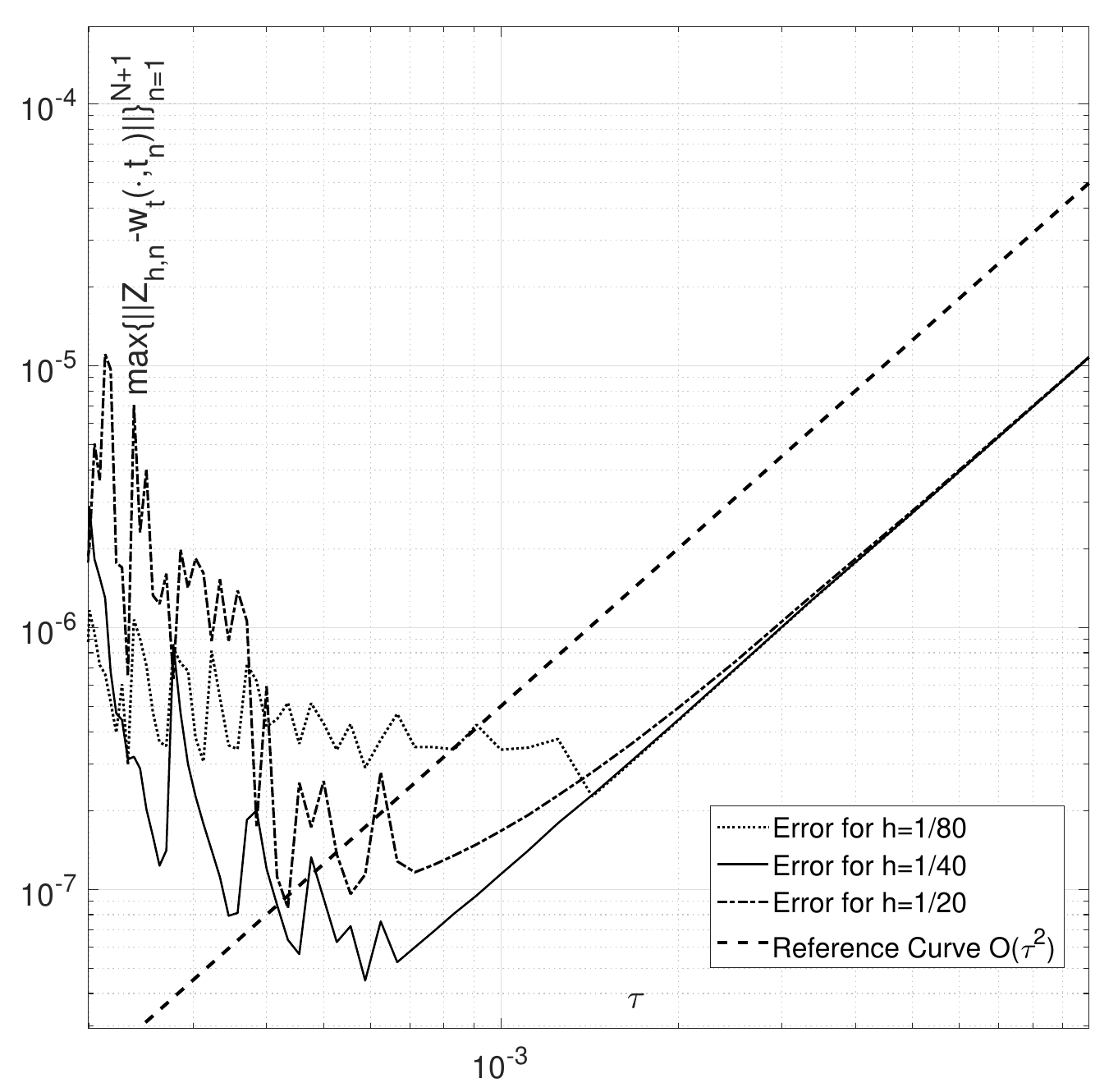}
	\caption{Error in log scale for $W_h$ (left) and $Z_h$ (right) throughout temporal step size $\tau$.}\label{fig_unsteady_err_tau}
\end{figure}

\begin{figure}[!htb]	 \centering \includegraphics[width=6.75cm,height=6.28cm]{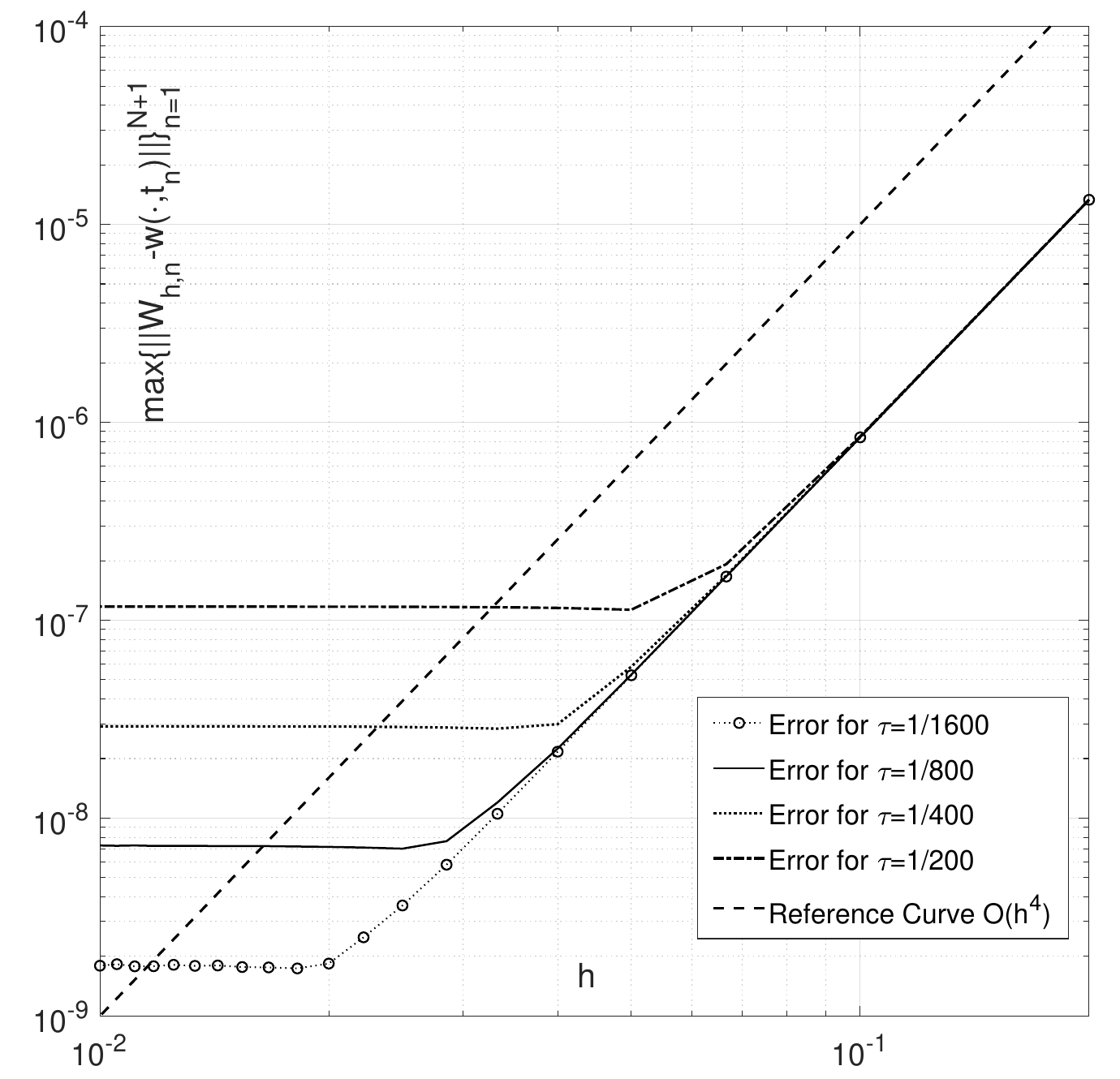}~~\includegraphics[width=6.75cm,height=6.25cm]{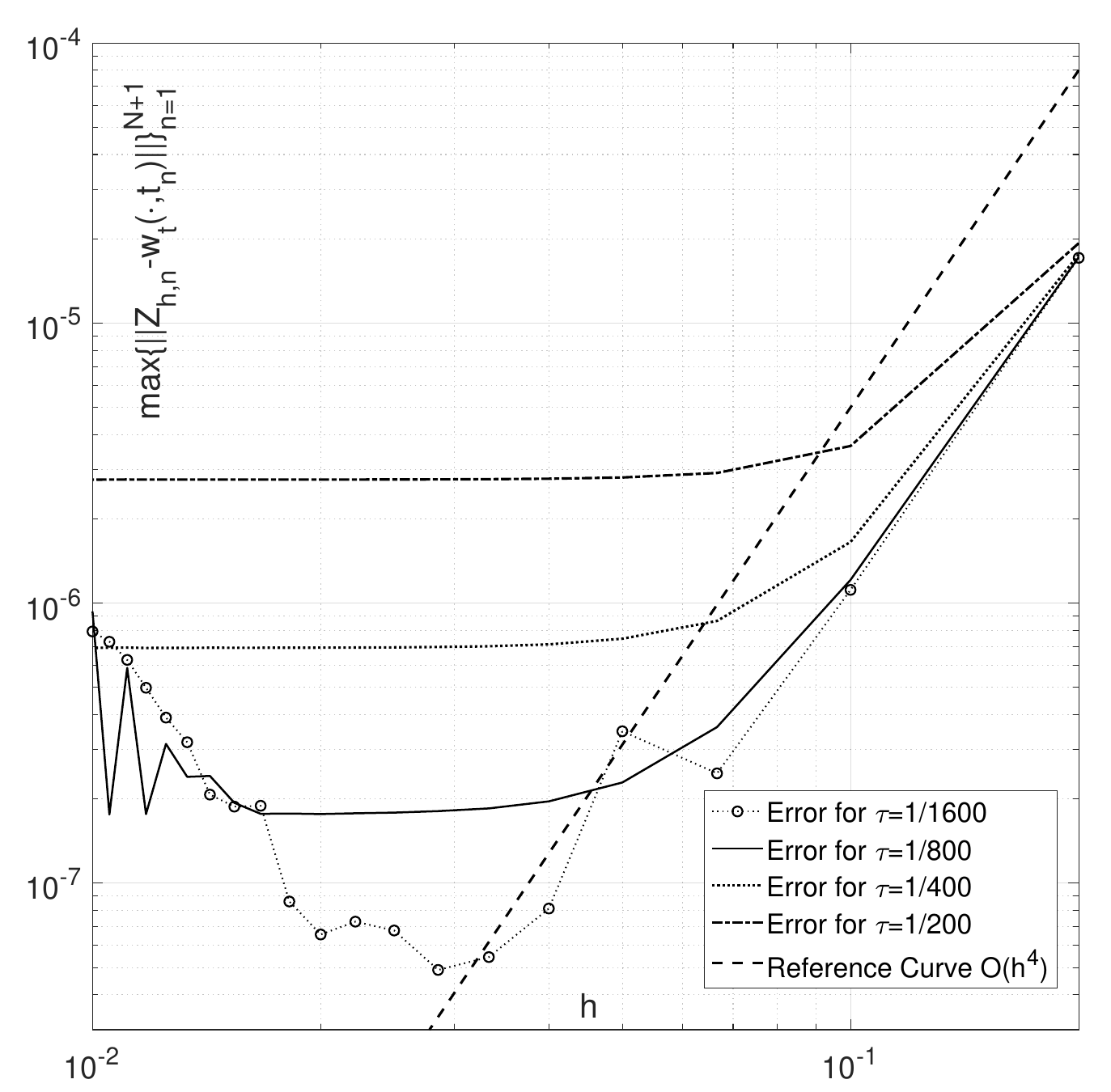}
	\caption{Error in log scale for $W_h$ (left) and $Z_h$ (right) throughout  spatial step size $h$.}\label{fig_unsteady_err_h}
\end{figure}

For hybrids methods, generally,  the order of accuracy cannot be determined directly by examining the components of the method independently. This requires a theoretical analysis, which is crucial yet often challenging. The resulting accuracy rates can be verified using specific test problems for which the exact solution is known. This allows for the optimization of the mesh parameters that support the theoretical findings. In Figures (\ref{fig_unsteady_err_tau}) and (\ref{fig_unsteady_err_h}), several step lengths are systematically utilized, and the relationship between these step lengths is examined to identify the optimal choice, which is found to be  $\tau\approx2~h^2$. Note that  $h\in[1/80,1/4]$ is suggested based on the numerical solution of the steady equation in the previous section. The behavior of the error on a logarithmic scale indicates an order of accuracy of $O(h^4)$ and $O(\tau^2)$ by choosing a fixed $\tau>0$ and $h>0$, respectively.

Following table confirms that the errors $W_h-w$ and $V_h-w_t$ in $C(0,T;L^2(0,l))$ norm  are $O(h^4+\tau^2)$. As noted in error plots, the error rate for  $V_h-w_t$ is deteriorated with the use of finer meshes.\\
\vskip .2in
\begin{tabular}{|c|c|c|c|c|}
  \hline
  $(h,\tau)$ & $(1/10,~1/50)$ &$(1/20,~1/200)$ & $(1/40,~1/800)$  & $(1/80,~1/3200)$ \\[3 pt] \hline 
   $||W_h-w||$ & $1.78\ 10^{-6}$ & $1.12\ 10^{-7}$ & $6.99\ 10^{-9}$ & $4.34\ 10^{-10}$\\ [3 pt]
  \hline
  $||V_h-w_t||$ & $4.48\ 10^{-5}$ & $2.79\ 10^{-6}$ & $1.78\ 10^{-7}$ & $3.07\ 10^{-7}$\\ [3 pt]
  \hline
\end{tabular}

\section{Conclusion}
In conclusion, the method proposed in this study can be seen as highly accurate and reliable.
 While several alternative methods exist in the literature and similar techniques have been employed for the wave equation, we believe that this computational study addresses a significant gap and serves as a robust alternative for the numerical solution of the beam equation. The advantages of the proposed method can be summarized as follows:

\begin{itemize}
  \item The discretization of the temporal derivative is performed independently, allowing for flexible implementation at this initial stage. Depending on the requirements, either a more accurate or more stable finite difference discretization can be employed.  Moreover, this flexibility in the construction of the method allows it to be adapted to more complicate problems such as non-linear beam  equation or plate equation. 
  \item The approximations of the derivatives $w_x$, $w_t$ and $w_{xt}$ are inherently derived from the algorithm. This feature eliminates the need for an additional discretization so it provides more stable calculations of synthetic measurement data related to derivative terms (velocity of deflection and slope).
  \item The logarithmic figures indicate that each component of the proposed hybrid method maintains its respective order of accuracy within the main algorithm. Notably, the optimization of the mesh parameters allows for the achievement of very low error values while utilizing a minimal number of mesh points in time and small number of elements in space."
\end{itemize}




\begin{thebibliography}{6}

\bibitem{VIN2007} Vinod K G, Gopalakrishnan S and  Ganguli R, {\it Free vibration and wave propagation analysis of uniform and tapered rotating beams using spectrally formulated finite elements.}  International Journal of Solids and Structures, \textbf{44} 5875-5893, 2007.
\bibitem{ANTOG:2000} Antognozzi M, {\it Investigation of the shear force contrast mechanism in transverse dynamic force microscopy} Ph.D. Thesis (UK: University of Bristol), 2000. 
\bibitem{NGU15} Nguyen T and et. al. {\it  Estimation of the shear force in transverse dynamic force microscopy using a sliding mode observer.}  AIP Advances, \textbf{5} 097157, 2015.
\bibitem{GB:XX2016} Bao G and Xu X, {\it An inverse random source problem in quantifying the elastic modulus of nanomaterials.} Inverse Problems, \textbf{29} 015006, 2013.
\bibitem{SG&DC&KJ&AS} Gunakala S R, Comissiong D M G, Jordan K. and Sankar A, {\it A Finite Element Solution of the Beam Equation via MATLAB.} International Journal of Applied Science and Technology, Vol. 2 No. 8, 2012.
\bibitem{WK&HB00}  Kwon Y W and Bang H, {\it The Finite Element Method Using Matlab.}  CRC Press, 2000.
\bibitem{OB&AH&AK24} Hasanov A, Kawano A and Baysal O, {\it Exponential stability of damped Euler-Bernoulli beam controlled by boundary springs and dampers},  Journal of Mathematical Analysis and Applications, Volume 533, Issue 2, 128031, 2024.
\bibitem{AH&OB16} Hasanov A and  Baysal O, {\it Identification of unknown temporal and spatial load distributions in a vibrating Euler-Bernoulli beam from Dirichlet boundary measured data.} Automatica 71 (2016) 106-117, 2016.
\bibitem{AH&OB&HI:19} Hasanov A,  Baysal O and Itou H, {\it Identification of an unknown shear force in a cantilever Euler-Bernoulli beam from measured boundary bending moment} J. Inverse Ill-posed Probl.  27(6), 859-876, 2019.
\bibitem{AH&AK&OB24}
Hasanov A, Kawano A and  Baysal O, {\it  Reconstruction of shear force in Atomic Force Microscopy from measured displacement of the cone-shaped cantilever tip.} Mathematics in Engineering, 6(1), 137–154, 2024.

\bibitem{KSSJNR17} Surana K S and Reddy J N {\it The Finite Element Method for Initial Value Problems
Mathematics and Computations.} CRC Press, 2017.

\bibitem{BH11} Brezis H, {\it Functional Analysis, Sobolev Spaces and Partial Differential Equations.} Springer, 2011.
\bibitem{ES&DM12} Süli E and Mayers D F, {\it An Introduction to Numerical Analysis.} Cambridge University Press, 2012.
\bibitem{AM12} Mayer A, {\it A Simplified Calculation of Reduced HCT–Basis Functions in a Finite Element Context}, Computational Methods in Applied Mathematics 12-4, 486-499, 2012.
\bibitem{MA20} Asadzadeh M, {\it An Introduction to the Finite Element Method (FEM) for Differential Equations.} Wiley, 2020.
\bibitem{OB&AH19} Baysal O and Hasanov A, {\it Solvability of the clamped Euler-Bernoulli beam equation.} Applied Mathematics Letters, vol. 93, 85-90, 2019.
\bibitem{AH&AGR2021} Hasanov A and Romanov A G, {\it Introduction to Inverse Problems for Differential Equations.} 2nd ed, Springer, New York, 2021.
\bibitem{ELC02} Evans L C, {\it Partial Differential Equations.} Rhode Island: American Mathematical Society, 2002.

%
%
%
%
%
%
%
%
%
%
%
%
%
%
%
%
%
%
%
%
%
%
%
%
%
%
%
%
%
%
%
%
%
%
%
%
%
%
%
%
%
%
%
\end{thebibliography}
\end{document}